\documentclass[11pt]{amsart}

\usepackage{amssymb,amsmath,amsthm, amsfonts, enumerate}
\newtheorem{thm}{Theorem}[section]
\newtheorem{prop}[thm]{Proposition}
\newtheorem{lem}[thm]{Lemma}
\newtheorem{cor}[thm]{Corollary}
\theoremstyle{definition}

\theoremstyle{remark}

\numberwithin{equation}{section}

\begin{document}

\title[Bound States for $(N+1)$-coupled LW-SW Equations]{
Existence of bound states for $(N+1)$-coupled long-wave--short-wave
interaction equations}

\author[Sharad Silwal]{Sharad Silwal}
\address{Jefferson College of Health Sciences,
101 Elm Ave SE, Roanoke, VA 24013, USA} \email{sdsilwal@jchs.edu}

\begin{abstract}
We prove the existence of an infinite family of smooth positive bound states
for $(N+1)$-coupled long-wave--short-wave interaction equations.
The system describes the interaction between $N$ short waves and a long wave and is of interest in physics and fluid dynamics.
\end{abstract}

\thanks{\textit{Date}: \today.}

 \maketitle

\section{Introduction}

%
%
%

\noindent A non-linear system of interaction between a complex
short-wave field $u$ and a real long-wave field $v$ of the form
\begin{equation}\label{singlelwsw}
\left\{
\begin{aligned}
  & i\partial_t u + \partial_{x}^{2} u = \alpha uv+\beta |u|^2u \\
   & \partial_t v+ \partial_x^{3}v+ v \partial_x v = \gamma \partial_{x} \left(|u|^2\right),\\
\end{aligned}
\right.
\end{equation}
was first studied in \cite{[Tsutsu]} concerning the well-posedness of the Cauchy problem. The system \eqref{singlelwsw}
which has an interaction between a nonlinear Schr\"{o}dinger
(NLS)-type short wave and a Korteweg de Vries (KdV)-type long wave
appears in a wide variety of physical systems. The reader may refer
to \cite{[Ben]} for a general theory of the non-linear
long-wave--short-wave interaction (LSI) model. Numerous successful
applications of the LSI model exist in different contexts of fluid
dynamics such as capillary-gravity waves in \cite{[KSK]},
sonic-Langmuir solitons in \cite{[Kar], [YO]}, Alfv\'{e}n waves in
\cite{[Sul]}, and Bose-Einstein condensates in \cite{[Nista]}, to
mention but a few.

Let us consider the multicomponent LSI system
\begin{equation}\label{lwsw}
\left\{
\begin{aligned}
  & i\partial_t u_j + \partial_{x}^{2} u_j = \alpha_j u_jv+\beta_j|u_j|^2u_j , \ j=1,2,. \ .\ .\ N,\\
   & \partial_t v+ \partial_x^{3}v+ v \partial_x v = \frac{1}{2} \ \partial_{x} \left(\sum_{l=1}^{N}\alpha_l|u_l|^2\right),
\end{aligned}
\right.
\end{equation}
where $u_1,...,u_N$ are complex-valued functions of $(x,t),$ $v$ is
a real-valued function of $(x,t),$ and $\alpha_j, \beta_j$ are real
constants. Systems of the form \eqref{lwsw} arise in the event of a
nonlinear interaction between $N$ short waves, namely, $u_1,. .
.,u_N$, and a long wave $v$. For more details, readers may see \cite{[AA], [Cra], [Kar], [KSK]}.
The non-linear multicomponent LSI model has not been studied as
extensively as its single component counterpart. As such, it has
been generating a significant amount of interest in recent years.
Some examples of where a multicomponent LSI system may arise are
bio-physics \cite{[Davy]}, ferromagnetism \cite{[My]}, water waves
\cite{[Cra]}, nonlinear optics \cite{[Ohta0]}, elastic solid
mechanics \cite{[Erb]}, and acoustics \cite{[Saz1]}.

In this paper, we are concerned with finding bound-state solutions
to \eqref{lwsw}. We look for solutions of the form
\begin{equation}\label{phipsi}
\left\{
\begin{aligned}
& u_j(x,t)=e^{i\omega t+ikx}\phi_j(x-ct),\ j=1,...,N,\\
& v(x,t)=\psi(x-ct),
\end{aligned}
\right.
\end{equation}
where $\psi$ and $\phi_{j}$ are real-valued functions and $c=2k.$
Inserting \eqref{phipsi} into \eqref{lwsw}, we find
$(\phi_1,\ldots,\phi_N,\psi)$ satisfies the following system of
ordinary differential equations:
\begin{equation}\label{phipsieqn}
\left\{
\begin{aligned}
& -\phi_j''+ \sigma\phi_j=-\beta_j\phi_j^3-\alpha_j\phi_j\psi,\ j=1,...,N,\\
&
-\psi''+c\psi=\frac{\psi^2}{2}-\frac{1}{2}\sum_{l=1}^{N}\alpha_l\phi_l^2,
\end{aligned}
\right.
\end{equation}
where $\sigma=k^2+\omega$, and the primes denote derivatives with
respect to the variable $\xi=x-ct$.

\medskip

To state our main result, let us denote by $Y$ the product space
\begin{equation*}
Y=H^1\times H^1\times . . . \times H^1\ (\textrm{taken}\ {N+1}\
\textrm{times}),
\end{equation*}
and $\Delta=(u_1,...,u_N,v).$ For $d>0$ and $\lambda\geq 0,$ we set
\begin{equation}
X_{d, \lambda} = \{ \Delta\in Y:
d\|v\|_{L^2}^2+\sum_{j=1}^{N}\|u_j\|_{L^2}^2=\lambda\}
\end{equation}
and consider the minimization problem
\begin{equation}\label{minimiz}
I(\lambda)=\inf\{E(\Delta):\Delta\in X_{d, \lambda}\},
\end{equation}
where
\begin{equation}
E(\Delta)=\int_{\mathbb{R}}\left(
\sum_{j=1}^{N}\left((u_j^\prime)^2+\frac{\beta_j}{2}u_j^{4}+\alpha_jv
u_j^{2} \right)+(v^\prime)^2-\frac{1}{3}v^3\right) \ dx
\end{equation}
The following is our main result:
\begin{thm}\label{familyofsol}
Suppose that $\alpha_j, \beta_j<0$ for all $j=1,...,N,$ and let
\begin{equation}\label{MLdef}
\mathcal{M}_\lambda = \{\Delta \in
X_{d,\lambda}:I(\lambda)=E(\Delta)\}.
\end{equation}
Then the following results hold:

\smallskip

(i) For every $\lambda > 0,$ we have $\mathcal{M}_\lambda \neq
\emptyset.$

\smallskip

(ii) If $(\phi_1,...,\phi_N,\psi)\in \mathcal{M}_\lambda,$ then
$\psi \geq 0$ and $\phi_{j}\geq 0 $where at least one $\phi_{j}$ and $\psi$
are nontrivial. Furthermore, $\psi, \phi_j\in
H^\infty$ and decay exponentially at infinity.

\smallskip

(iii) There exists $\lambda^\ast>0$ such that if
$(\phi_1,...,\phi_N,\psi)\in \mathcal{M}_\lambda$ for
$\lambda>\lambda^\ast,$ then $\psi(x) >0$ and $\phi_{j}(x)>0$ for all $x\in\mathbb{R}$ for at least one
$j\in\{1,...,N\}.$

(iv) There exists a family
\begin{equation}
\left\{
\begin{aligned}
& u_{j,n}(x,t)=e^{i\omega_{n}t+ik_{n}x}\phi_{j,n}(x-c_nt),\ j=1,...,N,\\
& v_n(x,t)=\psi_n(x-c_nt)
\end{aligned}
\right.
\end{equation}
of non-trivial bound-state solutions to \eqref{lwsw} with
\begin{equation}
\lim_{n\to \infty}c_n=+\infty.
\end{equation}
Moreover, each $\psi_n$ and $\phi_{j,n}$ can be chosen as in parts
(ii) and (iii).
\end{thm}

Theorem~\ref{familyofsol} is a generalization of the analogous result proved previously in
\cite{[Dias], [Dias2]} for (1+1)-coupled system given by
\eqref{singlelwsw}. We also mention the papers \cite{[C1], [C2],[SB1], [SB3]} where
different techniques have been used to prove the existence of solutions
for coupled systems of long and short-wave equations. In \cite{[AB11]}, the
existence of a two-parameter family of disjoint sets of bound states has been proved
in the case when $N=1.$ We use similar techniques as the ones used
in \cite{[Dias], [Dias2]}, which uses the concentration compactness
argument of P.L. Lions \cite{[L1]} and relies on the results of
Berestycki and Lions \cite{[Ber]}, along with an estimate for the
Lagrange multiplier of the associated problem.

\smallskip

We remark at this point that the multiconstraint variational problem
\begin{equation}\label{varremark}
\inf\{E(\Delta):\Delta\in Y, \|u_j\|_{L^2}^2=a_j,\ j=1,...,N,
\|v\|_{L^2}^2=a\}
\end{equation}
can also be used to establish the existence of bound-state solutions. The
solution set of the problem \eqref{varremark} corresponds not only
to $(N+1)$-parameter bound-state solutions of \eqref{lwsw}, but also
provides a disjoint family of such bound states for each $a_j>0$ and
$a>0.$ However, application of the concentration compactness
technique in this situation is complicated by the fact that one
needs to establish the strict subadditivity property for the
$(N+1)$-parameter problem \eqref{varremark} in order to establish
the relative compactness of the minimizing sequence.
In \cite{[SB11]}, this problem has been solved under three constraints for the energy functional
associated with 3-coupled Schr\"{o}dinger systems.
To our knowledge, the compactness of minimizing sequence for such multiconstraint problems with more than three constraints remains an open question.

\section{Existence of Minimizers}\label{varsolution}

\noindent We first establish a few lemmas and propositions that will
be used in the proof of our main result.

\begin{lem}\label{leftbd}
For all $\lambda>0,$ one has $I(\lambda)>-\infty$.
\end{lem}
\noindent\textbf{Proof.} Let $\Delta\in X_{d, \lambda}$. Then, upon
using the Gagliardo-Nirenberg inequality, we have, for each $j$,
\begin{equation*}
\begin{aligned}
\int u_j^4\leq \|u_j\|_{L^4}^4\leq
C_1\|u_j'\|_{L^2}\|u_j\|_{L^2}^3\leq C_1\lambda^{3/2}\|u_j'\|_{L^2}
\end{aligned}
\end{equation*}
and
\begin{equation*}
\begin{aligned}
v^3\leq \|v\|_{L^3}^3\leq C_2\|v'\|_{L^2}^{1/2}\|v\|_{L^2}^{5/2}\leq
C_2\left(\frac{\lambda}{d}\right)^{5/4}\|v'\|_{L^2}^{1/2}.
\end{aligned}
\end{equation*}
Also, for each $j$,
\begin{equation*}
\begin{aligned}
vu_j^2\leq \frac{1}{2}v^2+\frac{1}{2}u_j^4\leq
\frac{1}{2}\|v\|_{L^2}^{2}+\frac{1}{2}\|u_j\|_{L^4}^4\leq
\frac{\lambda}{2d}+\frac{C_1}{2}\lambda^{3/2}\|u_j'\|_{L^2}.
\end{aligned}
\end{equation*}
Using all of the above estimates, we finally obtain
\begin{equation*}
\begin{aligned}
E(\Delta)&=\sum_{j=1}^N\left(\|u_j'\|_{L^2}^2+\frac{\beta_j}{2}\int u_j^4+\alpha_j\int v u_j^{2}\right)+\|v'\|_{L^2}^2-\frac{1}{3}\int v^3\\
 &\geq\sum_{j=1}^N\left(\|u_j'\|_{L^2}^2-\frac{|\beta_j|}{2}\int u_j^4-|\alpha_j|\int |v| u_j^{2}\right)+\|v'\|_{L^2}^2-\frac{1}{3}\int v^3\\
 &\geq\sum_{j=1}^N\left(\|u_j'\|_{L^2}^2-\left(\frac{|\beta_j|}{2}C_1\lambda^{3/2}+\left(|\alpha_j|\left(\frac{\lambda}{2}+\frac{C_1}{2}\lambda^{3/2}\right)\right)\|u_j'\|_{L^2}\right)\right)\\
 &\quad+\|v'\|_{L^2}^2-\frac{C_2}{3}\left(\frac{\lambda}{d}\right)^{5/4}\|v'\|_{L^2}^{1/2}.
\end{aligned}
\end{equation*}
This proves that $E(\Delta)$ is bounded below by an expression
depending only on $\lambda$ and $d$.\qed

\begin{prop}\label{neg}
For all $\lambda\geq 0, I(\lambda)\leq 0$. Also, there exists
$\lambda^*>0$ such that for all $\lambda>\lambda^*, I(\lambda)\leq
-A\lambda^2$, where $A$ is a positive constant independent of
$\lambda$.
\end{prop}
\noindent\textbf{Proof.} Let $\Delta_1=(u_1,0,\ldots,0)\in X_{d,
\lambda}$ where $u_1\in H^1$ is such that $\|u_1\|_{L^2}^2=\lambda$.
Then, since $\beta_1<0$, we have
\begin{equation*}
\begin{aligned}
E(\Delta_1)=\int\left( u_1'^2+\frac{\beta_1}{2}u_1^4\right)\leq \int
u_1'^2.
\end{aligned}
\end{equation*}
Taking infimum over all $u_1\in H^1$ such that $\|u_1\|_{L^2}^2=\lambda$, we get $I(\lambda)\leq 0$.\\
Next, let $u_1\in H^1$ such that $\|u_1\|_{L^2}^2=1$ and put
$u_{1\lambda}=\lambda^{1/2}u_1$. Then,
$\Delta_{1\lambda}=(u_{1\lambda},0,\ldots,0)\in X_{d, \lambda}$.
Also, using the fact that $\beta_1<0$, we can have
\begin{equation*}
\begin{aligned}
E(\Delta_{1\lambda})=\lambda\int
u_1'^2+\frac{\beta_1}{2}\lambda^2\int u_1^4=\lambda\int
u_1'^2-\frac{|\beta_1|}{4}\lambda^2\int
u_1^4-\frac{|\beta_1|}{4}\lambda^2\int u_1^4.
\end{aligned}
\end{equation*}
Then, choosing $A=\frac{|\beta_1|}{4}\int u_1^4$ and $\lambda^*$ in
such a way that $\int u_1'^2-A\lambda^*\leq 0$, we will have, for
all $\lambda>\lambda^*$,
\begin{equation*}
\begin{aligned}
E(\Delta_{1\lambda})=\lambda\int u_1'^2-A\lambda^2-A\lambda^2\leq
\lambda A(\lambda^*-\lambda)-A\lambda^2\leq -A\lambda^2.
\end{aligned}
\end{equation*}
Taking infimum over all $u_1\in H^1$ such that $\|u_1\|_{L^2}^2=1$,
we get the result. \qed

\begin{lem}\label{positivesequence}
For all $\Delta=(f_1,...,f_N,g)\in Y,$ we have
\begin{equation*}
E(|\Delta|)=E(|f_1|,...,|f_N|,|g|)\leq E(f_1,...,f_N,g)=E(\Delta).
\end{equation*}
\end{lem}
\noindent\textbf{Proof.}
It is immediate from the following two facts: one is that, for $f\in
H^1$ real-valued, we have $\||f|'\|_{L^2}\leq \|f'\|_{L^2}$ (for example, see \cite{[Ber]}), and the
other that $\alpha_j<0$ yields $\alpha_j u_j^2|v|\leq \alpha_j u_j^2
v$ and $-\frac{1}{3}|v|^3\leq -\frac{1}{3}v^3$.\qed

\begin{lem}\label{allnegcutoff}
$I$ is non-increasing on $[0,\infty)$. That is, for all $\theta>1$, we have
$I(\theta\lambda)\leq\theta I(\lambda)$.
\end{lem}
\noindent\textbf{Proof.} Consider a sequence $\{\Delta_k\}=\{(u_{1,k},\ldots,u_{N,k},v_k)\}\subset X_{d, \lambda}$ and denote
$\sqrt{\theta}\Delta_k=\left(\sqrt{\theta}u_{1,k},\ldots,\sqrt{\theta}u_{N,k},\sqrt{\theta}v_k\right)$. Then
 \begin{equation*}
\begin{aligned}
&\quad \, E(\sqrt{\theta}\Delta_k)\\
&=\theta E(\Delta_k)-\frac{(\theta^{3/2}-\theta)}{3}\int v_k^3\\
&\qquad\qquad\qquad +\sum_{j=1}^N\left((\theta^{3/2}-\theta)\alpha_j\int v_k u_{j,k}^{2}+\frac{(\theta^{2}-\theta)}{2}\beta_j\int u_{j,k}^4\right)\\
&\leq \theta E(\Delta_k)\\
&-\min\{(\theta^{3/2}-\theta),(\theta^{2}-\theta)\}\left[\frac{1}{3}\int v_k^3+\sum_{j=1}^N\left(\frac{|\beta_j|}{2}\int u_{j,k}^4+|\alpha_j|\int v_k u_{j,k}^{2}\right)\right]\\
&\leq \theta E(\Delta_k).
\end{aligned}
\end{equation*}
Next, taking infimum over all sequences
$\Delta_k$ and using the fact that $\Delta\in X_{d, \lambda}$ if and
only if $\sqrt{\theta}\Delta\in X_{d,\theta\lambda}$, we obtain
$I(\theta\lambda)\leq\theta I(\lambda)$. \qed

\begin{lem}\label{nondecre}
There exists $\lambda_1\in[0,\lambda^*)$ such that $I(\lambda)<0 \Leftrightarrow \lambda>\lambda_1$. Furthermore, for all $\lambda>\lambda_1$ and $\theta>1$, we have
$I(\theta\lambda)<\theta I(\lambda)$.
\end{lem}
\noindent\textbf{Proof.} By Proposition \ref{neg}
and Lemma \ref{allnegcutoff}, $I$ is a non-positive and non-increasing function
on $[0,\infty)$ which is strictly negative on $(\lambda^*,\infty)$.
Hence, there exists $\lambda_1\in[0,\lambda^*]$ such that $I$ vanishes on $[0,\lambda_1]$ and is strictly negative and decreasing on $(\lambda_1,\infty)$.
Then, for any $\lambda>\lambda_1$, we have $I(\lambda)<0$. Hence, there exists
$\delta>0$ such that we can, as guaranteed by Lemma
\ref{positivesequence}, construct a positive
 minimizing sequence \[\Delta_k=(u_{1,k},\ldots,u_{N,k},v_k)\in X_{d, \lambda}\] for the minimization problem \eqref{minimiz} in such a way that $E(\Delta_k)\leq -\delta$ for all $k\in \mathbb{N}\setminus\{0\}$.
 It exists because otherwise it would give us a subsequence $\Delta_{k_l}$ with $E(\Delta_{k_l})\geq 0$ which, in turn, would lead to a contradiction
 $I(\lambda)=\displaystyle\lim_{k_l\rightarrow\infty} E(\Delta_{k_l})\geq 0$. So, using the fact that $v_k>0, \alpha_j<0$ and $\beta_j<0$, we get

 \begin{equation*}
\begin{aligned}
-\delta &\geq E(\Delta_k)\\
&=\sum_{j=1}^N\left(\|u_{j_k}'\|_{L^2}^2+\frac{\beta_j}{2}\int u_{j,k}^4+\alpha_j\int v_k u_{j,k}^{2}\right)+\|v_k'\|_{L^2}^2-\frac{1}{3}\int v_k^3\\
&\geq -\sum_{j=1}^N\left(\frac{|\beta_j|}{2}\int
u_{j,k}^4+|\alpha_j|\int v_k u_{j,k}^{2}\right)-\frac{1}{3}\int
v_k^3.
\end{aligned}
\end{equation*}
Combining this result with what we have computed in Lemma \ref{allnegcutoff}, we obtain
 \begin{equation*}
\begin{aligned}
&\quad  E(\sqrt{\theta}\Delta_k)\\
&\leq \theta E(\Delta_k)\\
&-\min\{(\theta^{3/2}-\theta),(\theta^{2}-\theta)\}\left[\frac{1}{3}\int v_k^3+\sum_{j=1}^N\left(\frac{|\beta_j|}{2}\int u_{j,k}^4+|\alpha_j|\int v_k u_{j,k}^{2}\right)\right]\\
&\leq \theta
E(\Delta_k)-\delta\min\{(\theta^{3/2}-\theta),(\theta^{2}-\theta)\}\\
&<\theta E(\Delta_k).
\end{aligned}
\end{equation*}
The result follows upon taking infimum over all positive minimizing sequences
$\Delta_k$ and using the fact that $\Delta\in X_{d, \lambda}$ if and
only if $\sqrt{\theta}\Delta\in X_{d,\theta\lambda}$. \qed

\medskip

With Lemma \ref{nondecre} in hand, an argument
very similar to Lemma 2.3 of \cite{[Ohta0]} yields the following strict sub-additivity result for $I(\lambda)$:
\begin{cor}\label{subadd}
Let $\lambda>\lambda_1$ and $0<\Omega<\lambda$. Then
$I(\lambda)<I(\Omega)+I(\lambda-\Omega)$.
\end{cor}

The following Proposition establishes the existence of minimizers
for the minimization problem \eqref{minimiz}:
\begin{prop}\label{minimizexist}
For every $\lambda>\lambda_1\geq 0,$ the set $\mathcal{M}_\lambda$
as defined in \eqref{MLdef} is non-empty. Furthermore, if
$(\phi_1,\ldots,\phi_N,\psi)\in \mathcal{M}_\lambda,$ then the
following hold:
\begin{enumerate}[(i)]
\item $\psi \geq 0$ and $\phi_{j}\geq 0 $ where $\psi$ and at least one $\phi_{j}$
are nontrivial.
\item $\psi, \phi_j\in H^\infty$ and decay exponentially at infinity.
\end{enumerate}
\end{prop}
\noindent\textbf{Proof.} Consider a positive
 minimizing sequence $\Delta_k=(u_{1,k},\ldots,u_{N,k},v_k)\in X_{d, \lambda}$ for the minimization problem \eqref{minimiz}. Set \[\rho_k=\sum_{j=1}^N u_{j,k}^2+dv_k^2\] and,
 following the notations in \cite{[L1]}, define the concentration function of $\rho_k$ as \[Q_k(r)=\sup_{y\in\mathbb{R}}\int_{y-r}^{y+r}\rho_k, \, \text{and set}\: \Omega(r)=\lim_{k\rightarrow\infty} Q_k(r)\:\text{and}\: \Omega=\lim_{r\rightarrow\infty} \Omega(r).\]
 Then there are three possibilities for $\Omega$: vanishing $(\Omega =0)$, dichotomy $(0<\Omega<\lambda)$ and compactness $(\Omega=\lambda)$. Our goal is to establish the last alternative. We now divide the proof into three parts.

 \smallskip

 \textbf{Part I. The vanishing case does not occur.} If $\Omega=0$, then, by the definition of $\rho_k$, we have, for each $j=1,\ldots, N$,
 \[\lim_{k\rightarrow\infty}\sup_{y\in\mathbb{R}}\int_{y-r}^{y+r}u_{j,k}^2=\lim_{k\rightarrow\infty}\sup_{y\in\mathbb{R}}\int_{y-r}^{y+r} v_{k}^2=0.\]
 Using Lemma I.1 of \cite{[L1]} and the fact that $u_{j,k}$ and $v_k$ are bounded in $H^1$ as seen in the proof of Lemma \ref{leftbd}, we get
 \[\|u_{j,k}\|_{L^p}\rightarrow 0, \|v_k\|_{L^p}\rightarrow 0,\: \text{for all}\: p>2.\]
 Then, upon using the Cauchy-Schwarz inequality, we have
 \[0\leq \lim_{k\rightarrow\infty}\int v_k u_{j,k}^{2}\leq \lim_{k\rightarrow\infty}\|v_k\|_{L^2}\|u_{j,k}\|_{L^4}^2=0.\]
 This implies that
\begin{equation}\label{finallimit}
\begin{aligned}
& \quad\: I(\lambda)
=\displaystyle\lim_{k\rightarrow\infty}\inf E(\Delta_k)\\
&=\lim_{k\rightarrow\infty}\inf\sum_{j=1}^N\left(\|u_{j,k}'\|_{L^2}^2-\frac{|\beta_j|}{2}\|u_{j,k}\|_{L^4}^4-|\alpha_j|\int v_k u_{j,k}^{2}\right)+\|v_k'\|_{L^2}^2-\frac{1}{3}\|v_k^3\|_{L^3}^3\\
&=\lim_{k\rightarrow\infty}\inf\sum_{j=1}^N\|u_{j,k}'\|_{L^2}^2+\|v_k'\|_{L^2}^2
\geq 0,
\end{aligned}
\end{equation}
which is contradictory to Proposition \ref{neg}, and the vanishing
case is ruled out.

\smallskip

\textbf{Part II. The dichotomy case does not occur.} Assume the
dichotomy case $(0<\Omega<\lambda)$. Then, using Lemma III.1 of
\cite{[L1]}, for all $\epsilon>0$, we can have, for all $k\in
\mathbb{N}$,
\[\left|\int \eta_1^2(x)\rho_k(x-y_k)-\Omega\right|<\epsilon \: \text{and}\: \left|\int \eta_2^2(x)\rho_k(x-y_k)-(\lambda-\Omega)\right|<\epsilon,\]
where a sequence $\{y_k\}\subset \mathbb{R}$ and cut-off functions
$\eta_1,\eta_2\in C^\infty (\mathbb{R})$ are so chosen that for some
fixed constants $0<R_1<\frac{R_2}{2}$,
\begin{enumerate}[(i)]
\item $0\leq \eta_i\leq 1, |\eta_i'|<\epsilon<\epsilon, i=1,2,$
\item $\eta_1(x)=1$ for $|x|\leq R_1$ and $\eta_1(x)=0$ for $|x|\geq \frac{R_2}{2}$,
\item $\eta_2(x)=1$ for $|x|\geq R_2$ and $\eta_2(x)=0$ for $|x|\leq \frac{R_2}{2}$,
\item $\int_{R_1\leq |x-y_k|\leq R_2} \rho_k(x) \leq \epsilon$.
\end{enumerate}
Set $u_{j,k}^{(i)}=\eta_i(x-y_k)u_{j,k}$ and
$v_{k}^{(i)}=\eta_i(x-y_k)v_{k}, i=1,2$. Then, using the Sobolev and
Cauchy-Schwarz inequalities, the fact that
\[u_{j,k}^2(x)\leq \rho(x) \: \text{and} \: dv_{k}^2(x)\leq \rho(x),\]
and the above inequality (iv), we obtain
\begin{equation*}
\begin{aligned}
&\quad\: \int_{R_1\leq |x-y_k|\leq \frac{R_2}{2}} \left[(u_{j,k}')^2-(u_{j,k}^{(1)'})^2\right]\\
&= \int_{R_1\leq |x-y_k|\leq \frac{R_2}{2}}\bigg[(1-\eta_1^2(x-y_k))(u_{j,k}')^2\\
&\qquad\qquad-\eta_1'^2(x-y_k)u_{j,k}^2-2\eta_1(x-y_k)\eta_1'(x-y_k)u_{j,k}u_{j,k}'\bigg]\\
&\geq -\epsilon^3-2C_1\epsilon=-\epsilon^3-C\epsilon.
\end{aligned}
\end{equation*}
The exact same inequality holds when $u_{j,k}^{(2)}$ is used instead
of $u_{j,k}^{(1)}$ in the above computation. Hence, we have
\begin{equation*}
\begin{aligned}
&\quad\: \|u_{j,k}'\|_{L^2}^2-\|u_{j,k}^{(1)'}\|_{L^2}^2-\|u_{j,k}^{(2)'}\|_{L^2}^2\\
&= \int_{R_1\leq |x-y_k|\leq R_2}u_{j,k}'^2-\int_{R_1\leq |x-y_k|\leq \frac{R_2}{2}}(u_{j,k}^{(1)'})^2-\int_{\frac{R_2}{2}\leq |x-y_k|\leq R_2}(u_{j,k}^{(2)'})^2\\
&=\int_{R_1\leq |x-y_k|\leq \frac{R_2}{2}} (u_{j,k}')^2-(u_{j,k}^{(1)'})^2+\int_{\frac{R_2}{2}\leq |x-y_k|\leq R_2} (u_{j,k}')^2-(u_{j,k}^{(2)'})^2\\
&\geq -2\epsilon^3-2C\epsilon.
\end{aligned}
\end{equation*}
The exact same computation yields
\[\|v_{k}'\|_{L^2}^2-\|v_{k}^{(1)'}\|_{L^2}^2-\|v_{k}^{(2)'}\|_{L^2}^2\geq-2\epsilon^3-2C\epsilon.\]
Furthermore, using the Gagliardo-Nirenberg inequality, we obtain
\[\|u_{j,k}\|_{L^4}^4-\|u_{j,k}^{(1)}\|_{L^4}^4-\|u_{j,k}^{(2)}\|_{L^4}^4\leq C\epsilon^3,\;
\|v_{k}\|_{L^3}^3-\|v_{k}^{(1)}\|_{L^3}^3-\|v_{k}^{(2)}\|_{L^3}^3\leq
C\epsilon^{5/2}.\] Next, Cauchy-Schwarz and Gagliardo-Nirenberg inequalities yield
\begin{equation*}
\begin{aligned}
&\quad\: \int v_k u_{j,k}^{2}-\int v_k^{(1)} (u_{j,k}^{(1)})^{2}-\int v_k^{(2)} (u_{j,k}^{(2)})^{2}\\
& =\int_{R_1\leq |x-y_k|\leq \frac{R_2}{2}}(1-\eta_1^3(x-y_k))v_k u_{j,k}^2\\
&\qquad\qquad+\int_{\frac{R_2}{2}\leq |x-y_k|\leq R_2}(1-\eta_2^3(x-y_k))v_k u_{j,k}^2\\
& \leq\int_{R_1\leq |x-y_k|\leq R_2}v_k u_{j,k}^2\leq
\|v_k\|_{L^2}\|u_{j,k}\|_{L^4}^2\leq C\epsilon^{5/4}.
\end{aligned}
\end{equation*}
Denoting
$\Delta_k^{(i)}=(u_{1,k}^{(i)},\ldots,u_{N,k}^{(i)},v_k^{(i)}),
i=1,2$, and putting together the above computations, we finally
obtain
\begin{equation*}
\begin{aligned}
&\quad\: E(\Delta_k)\\
& =\sum_{j=1}^N\left(\|u_{j,k}'\|_{L^2}^2-\frac{|\beta_j|}{2}\|u_{j,k}\|_{L^4}^4-|\alpha_j|\int v_k u_{j,k}^{2}\right)+\|v_k'\|_{L^2}^2-\frac{1}{3}\|v_k^3\|_{L^3}^3\\
&\geq E(\Delta_k^{(1)})+E(\Delta_k^{(2)})-C(\epsilon).
\end{aligned}
\end{equation*}
Hence, $I(\lambda)\geq I(\Omega)+I(\lambda-\Omega)$, which
contradicts Corollary \ref{subadd}.

\smallskip

\textbf{Part III. The set $\mathcal{M}_\lambda$ is non-empty.} Parts
I and II imply that the only remaining case is the last one,
namely, $\Omega=\lambda$. This implies that the sequence $\Delta_k$
is relatively compact up to translations. That is, there exist its
subsequence, again denoted $\Delta_k=(u_{1,k},\ldots,u_{N,k}, v_k)$,
a sequence $\{y_k\}\subset\mathbb{R}$ and the limit
$(\phi_1,\ldots,\phi_N,\psi)\in X_{d, \lambda}$ such that
\[\tilde{\Delta}_k=(\tilde{u}_{1,k},\ldots,\tilde{u}_{N,k}, \tilde{v}_k)\rightharpoonup(\phi_1,\ldots,\phi_N,\psi)\: \text{in}\: H^1(\mathbb{R}),\]
where $\tilde{u}_{j,k}=u_{j,k}(\cdot-y_k), j=1,\ldots, N$ and
$\tilde{v}_k=v_k(\cdot-y_k)$. Consequently, we have, for all $p\geq 2$,
\begin{equation}\label{weakinH1}
\tilde{\Delta}_k\rightarrow (\phi_1,\ldots,\phi_N,\psi)\: \text{in}\: L^p(\mathbb{R})
\: \text{ and}\:\int v_ku_{j,k}^2\rightarrow \int\psi\phi^2.
\end{equation}
Computing as in \eqref{finallimit} and using \eqref{weakinH1}, we obtain
\begin{equation}\label{stronginH1}
I(\lambda)=\liminf\limits_{k\rightarrow\infty}E(\tilde{\Delta}_k)=E((\phi_1,\ldots,\phi_N,\psi)).
\end{equation}
But, on the other hand, we have,
\[E((\phi_1,\ldots,\phi_N,\psi))\geq \inf_{\Delta\in X_{d, \lambda}}  E(\Delta)=I(\lambda).\]
This results in $E((\phi_1,\ldots,\phi_N,\psi))=I(\lambda)$ and we thus have
$(\phi_1,\ldots,\phi_N,\psi)\in \mathcal{M}_\lambda\neq\emptyset$.
The computation in \eqref{stronginH1} also establishes
\[\liminf\limits_{k\rightarrow\infty}\sum_{j=1}^N\|\tilde{u}_{j,k}'\|_{L^2}^2+\|\tilde{v}_k'\|_{L^2}^2
=\sum_{j=1}^N\|\phi_{j}'\|_{L^2}^2+\|\psi'\|_{L^2}^2,\] which leads
to
\[\tilde{\Delta}_k=(\tilde{u}_{1,k},\ldots,\tilde{u}_{N,k}, \tilde{v}_k)\rightarrow(\phi_1,\ldots,\phi_N,\psi)\: \text{in}\: H^1(\mathbb{R}).\]
Then, since $\tilde{u}_{j,k}, \tilde{v}_{k}\geq 0$, we must have
$\phi_j, \psi\geq 0$. Note that $(\phi_1,\ldots,\phi_N,\psi)\in
X_{d, \lambda}$ if and only if $(\phi_1,\ldots,\phi_N,\psi+a)\in
X_{d, \lambda}$, for any $a>0$. Furthermore, due to the terms
associated to $\psi$ being negative, we have
\[E((\phi_1,\ldots,\phi_N,\psi+a))< E((\phi_1,\ldots,\phi_N,\psi)),\] which, if $\psi= 0$, would be contradictory to the fact that \[E((\phi_1,\ldots,\phi_N,\psi))=I(\lambda)=\inf\{E(\Delta):\Delta\in X_{d, \lambda}\},\]
and hence, $\psi\not\equiv 0$ on $\mathbb{R}.$
Using the same argument as in the proof of Proposition 4.5 of \cite{[Dias]}, one can see that $\phi_j\not\equiv 0$
for at least one $j\in\{1,...,N\}.$
The smoothness of $\phi_{j}, j=1,\ldots,N$ and $\psi$ can be
established by a standard bootstrap technique. Next, choosing
$c=2k>0$ and $\sigma=k^2+\omega>0$, and using Theorem 8.1.1 in
\cite{[Caz]}, there exists a $\varepsilon>0$ such that
$e^{\varepsilon|\cdot|}\phi_{j}, e^{\varepsilon|\cdot|}\psi\in L^{\infty}$.
Hence the functions $\phi_{j}, j=1,\ldots,N$ and $\psi$ are in
$H^\infty$ and decay exponentially at infinity.\qed

\smallskip

\section{Proof of the main theorem}

\noindent In this section, we establish our main result. Parts (i)
and (ii) have already been established. To prove the remaining
statements, we first prove two ancillary results in the form of a
lemma and a proposition.

\begin{lem}\label{LMneg}
There exists a constant $A>0$ and $\lambda^\ast>0$ such that for all
$\lambda>\lambda^\ast,$ the Lagrange multiplier $\mu$ satisfies
\begin{equation*}
\mu \leq -A\lambda.
\end{equation*}
\end{lem}
\noindent\textbf{Proof.}
By Proposition~\ref{minimizexist}, there exists
$(\phi_1,\ldots,\phi_N,\psi)\in \mathcal{M}_\lambda$ such that
$\psi>0$ and $\phi_j\geq 0, j=1,\ldots, N$. There exists a Lagrange
multiplier $\mu$ depending only on $\lambda$ such that
\begin{equation*}
\left\{
\begin{aligned}
& -\phi_j''-\mu\phi_j=-\beta_j\phi_j^3-\alpha_j\phi_j\psi,\ j=1,...,N,\\
&
-\psi''-\mu d\psi=\frac{\psi^2}{2}-\frac{1}{2}\sum_{l=1}^{N}\alpha_l\phi_l^2.
\end{aligned}
\right.
\end{equation*}
Note that there are $N+1$ equations in total, keeping in mind that
one or more of them, determined by how many of $\phi_j, j=1,\ldots,
N$ are trivial, might vanish. Multiplying each of the $N+1$
equations by $\phi_j, j=1,\ldots, N$ and $\psi$ respectively, and
integrating by part,
\begin{equation*}
\left\{
\begin{aligned}
& \int \phi_j'^2-\mu\int\phi_j^2=-\beta_j\int\phi_j^4-\alpha_j\int\phi_j^2\psi,\ j=1,...,N,\\
& \int \psi'^2-\mu d\int\psi^2=\frac{1}{2}\int
\psi^3-\frac{1}{2}\int\psi\sum_{l=1}^{N}\alpha_l\phi_l^2.
\end{aligned}
\right.
\end{equation*}
Adding,
\begin{equation*}
\begin{aligned}
&\quad\sum_{j=1}^{N}\|\phi_j'\|_{L^2}^2+\|\psi'\|_{L^2}^2-\mu\lambda\\
&=\frac{1}{2}\|\psi\|_{L^3}^3-\int\left[\sum_{j=1}^{N}\left(\alpha_j\phi_j^2\psi+\beta_j
\phi_j^4+\frac{\psi}{2}\sum_{l=1}^{N}\alpha_l\phi_l^2\right)\right].
\end{aligned}
\end{equation*}
On the other hand, we have
\begin{equation*}
\begin{aligned}
&\quad I(\lambda)=E((\phi_1,\ldots,\phi_N,\psi))=\\
&\sum_{j=1}^N\left(\|\phi_{j}'\|_{L^2}^2+\frac{\beta_j}{2}\|\phi_{j}\|_{L^4}^4+\alpha_j\int
\psi
\phi_{j}^{2}\right)+\|\psi'\|_{L^2}^2-\frac{1}{3}\|\psi\|_{L^3}^3
\end{aligned}
\end{equation*}
Hence,
\begin{equation*}
\begin{aligned}
\quad \mu\lambda
& =I(\lambda)+\sum_{j=1}^N\frac{1}{2}\int\sum_{l=1}^{N}\alpha_l\phi_l^2\psi+\frac{\beta_j}{2}\|\phi_{j}\|_{L^4}^4-\frac{1}{6}\|\psi\|_{L^3}^3\\
&\leq I(\lambda)\leq -A\lambda^2,
\end{aligned}
\end{equation*}
for some $A>0$ independent of $\lambda$, and for all
$\lambda>\lambda^*$, as given by Proposition \ref{neg}.
\qed

\begin{prop}\label{allpos}
There exists $\lambda^\ast>0$ such that if
$(\phi_1,...\phi_N,\psi)\in \mathcal{M}_\lambda$ for
$\lambda>\lambda^\ast,$ then we have $\psi(x)>0$ and $\phi_j(x)>0$ for all $x\in\mathbb{R}$ for at least one $j\in\{1,...,N\}.$
\end{prop}
\noindent\textbf{Proof.}
Let $A$ and $\lambda^\ast$ be as in Lemma~\ref{LMneg}. Then
$s=-\mu>0.$ For each $j=1,...,N,$ the functions $\psi$ and $\phi_j$ satisfy
\begin{equation}\label{phipos}
\left\{
\begin{aligned}
& -\phi_j^{\prime \prime}+s\phi_j=-\beta_j\phi_j^3-\alpha_j\phi_j\psi,\ j=1,2,...,N,\\
  & -\psi''+sd\psi=\frac{\psi^2}{2}-\frac{1}{2}\sum_{l=1}^{N}\alpha_l\phi_l^2.
   \end{aligned}
\right.
\end{equation}
Let $P_s(x)=\frac{1}{2\sqrt{s}}e^{-\sqrt{s}|x|}$ for $x\in
\mathbb{R}.$ Then \eqref{phipos} can be written as
\begin{equation*}
\left\{
\begin{aligned}
& \phi_j=P_s\star \left( b_j\phi_j^3+a_j\phi_j \psi\right),\ j=1,2,...,N,\\
  & \psi=P_{sd}\star \left(\frac{\psi^2}{2}+\frac{1}{2}\sum_{l=1}^{N}a_l\phi_l^2\right).
   \end{aligned}
   \right.
\end{equation*}
where $b_j=-\beta_j>0$ and $a_j=-\alpha_j>0.$
Since the convolution
of the positive kernel $P_s$ with a nonnegative and not identically zero function always
produces a positive function, it follows that $\phi_j$ is
positive on all of $\mathbb{R}$ for at least one $j\in\{1,...,N\}.$
\qed
\medskip

We now complete the proof of our main result.
\medskip

\noindent\textbf{Proof of Theorem \ref{familyofsol}:} Parts (i) and
(ii) have been established in Proposition~\ref{minimizexist} while
part (iii) has been proved in Proposition~\ref{allpos}. Next, we
establish the remaining part (iv). By Lemma~\ref{LMneg}, we have
\[\mu\leq-A\lambda<0.\]
Now, choose a sequence $\lambda_n\rightarrow \infty$. For each $n$,
there will be a Lagrange multiplier $\mu_n$ (depending only on
$\lambda_n$) such that $\mu_n\rightarrow-\infty$. Corresponding to
each Lagrange multiplier $\mu_n$ are solutions $\phi_{j,n},
j=1,\ldots,N$ and $\psi_n$. Set \[c_n=-\mu_n, \:
k_n=-\frac{1}{2}\mu_n,\:\omega_n=-\mu_n-k_n^2,\:
v_n(x,t)=\psi_n(x-c_nt)\] and
\[u_{j,n}(x,t)=e^{i\omega_{j,n}t+ik_{n}x}\phi_{j,n}(x-c_nt),\ j=1,...,N.\]
Then, indeed,
\[\lim_{n\to \infty}c_n=\lim_{n\to \infty}-\mu_n=+\infty.\]
Furthermore, since $k_n=-\frac{1}{2}\mu_n>0$ and
$\sigma_n=k_n^2+\omega_n=-\mu_n>0$, a similar argument as in the
proof of Proposition~\ref{minimizexist} establishes that the
functions $\phi_{j,n}, j=1,\ldots,N$ and $\psi_n$ are in $H^\infty$
and decay exponentially at infinity.

\qed

\medskip


\section{Acknowledgments}
\noindent The author is  partially  supported  by a grant under JCHS
Faculty Development Grant Program.

\footnotesize


\end{document}